\documentclass[12pt]{amsart}
\usepackage{amsthm}
\usepackage{amsmath}
\usepackage{amsfonts} 

\title{\bf  Binomial Residues}
\author{Eduardo Cattani, Alicia Dickenstein and Bernd Sturmfels}
\address{Eduardo Cattani: Department of Mathematics and Statistics. 
University
of Massachusetts. Amherst, MA 01003, USA}
\email{cattani@math.umass.edu}

\address{Alicia Dickenstein: Departamento~de Matematica.
Universidad de Buenos Aires. (1428) Buenos Aires, Argentina}
\email{alidick@dm.uba.ar}

\address{Bernd Sturmfels: Department~of Mathematics. 
University of California.
Berkeley, CA 94720, USA}
\email{bernd@math.berkeley.edu}

\newcommand{\baseRing}[1]{\ensuremath{\mathbb{#1}}}
\newcommand{\Z}{\baseRing{Z}}
\newcommand{\R}{\baseRing{R}}
\newcommand{\C}{\baseRing{C}}
\newcommand{\N}{\baseRing{N}}
\newcommand{\Q}{\baseRing{Q}}

\newcommand{\CP}{\baseRing{P}}

\def\pd#1{ \partial_{#1} }

\theoremstyle{plain}
\newtheorem{theorem}{Theorem}[section]
\newtheorem{lemma}[theorem]{Lemma}
\newtheorem{corollary}[theorem]{Corollary}
\newtheorem{prop}[theorem]{Proposition}
\newtheorem{proposition}[theorem]{Proposition}

\theoremstyle{definition}
\newtheorem{algorithm}[theorem]{Algorithm}
\newtheorem{definition}[theorem]{Definition}

\numberwithin{equation}{section}

\DeclareMathOperator{\nsupp}{nsupp}

\newcommand{\rx}{{\rm Res}^X}

\newcommand{\La}{M}

\newcommand{\Script}[1]{\ensuremath{{\mathcal{#1}}}}
\newcommand{\HH}{\Script{H}}

\newcommand{\ScS}{\Script{S}}

\newcommand{\II}{\Script{I}}

\newcommand{\scA}{\Script{A}}

\newcommand{\one}{\mathbf{1}}

\def\proof{\smallskip\noindent{\sl Proof.\ }}
\begin{document}

\begin{abstract} 
A binomial residue is a rational function defined by a
hypergeometric integral whose kernel is singular along binomial
divisors. Binomial residues provide an integral representation
for  rational solutions of $A$-hypergeometric systems of
Lawrence type. The space of binomial residues of a given degree,
 modulo those which are polynomial in some variable, has
dimension equal to the Euler characteristic of the 
matroid associated with $A$.
\end{abstract}

\maketitle

\section{Introduction}

By a {\it binomial residue} we mean a rational  function 
in $2n$ variables $x_1,\dots \! ,x_n$, $ y_1,\dots\!,y_n $,
which is defined by a  residue integral of the form
\begin{equation}
\label{localresidue}
R_\Gamma(x,y) \ :=\ \int_\Gamma 
 \frac{t^{\gamma} }
{ (x_1 + t^{a_1} y_1)^{\beta_1}
\cdots (x_n + t^{a_n} y_n)^{\beta_n}} \frac{dt_1}{t_1} \wedge 
\dots \wedge
\frac{dt_d}{t_d}.
\end{equation}

\noindent Here $a_1,a_2,\dots,a_n$ are 
non-zero lattice vectors which 
span $\Z^d$,  $\gamma$ is any vector in $\Z^d$,
$\beta_1,\ldots,\beta_n$ are positive integers,
and $\Gamma$ ranges over a certain collection,
 specified in  (\ref{singleintegral}) below,
of compact $d$-cycles in the torus $(\C^*)^d$.
In this paper we study analytic, combinatorial, and
geometric properties of binomial residues.  On the
analytic side, we view binomial residues as hypergeometric 
integrals  \cite[page
223]{sst2} and, consequently, as rational solutions of a
certain 
$A$-hypergeometric system of differential equations, 
in the sense of
Gel'fand, Kapranov and Zelevinsky
\cite{gkz89, gkz90}.  The $A$-hypergeometric system annihilating
(\ref{localresidue}) is
the left  ideal in the $2n$-dimensional
 Weyl algebra  generated by the operators
\begin{eqnarray}\label{system}
& 
\partial_x^u\,\partial_y^v\, -\,\partial_x^v\,\partial_y^u
\,\,\, \hbox{whenever} \,\,\,u,v\in\N^n 
\,\,\, \hbox{and} \,\,\,\sum_{i=1}^n u_i a_i =  
\sum_{i=1}^n v_i a_i  ,\nonumber \\
& x_i \partial_{x_i} + y_i \partial_{y_i} + \beta_i
\quad \hbox{for} \quad i = 1,2,\ldots,n , \qquad
\, \hbox{and}\\ &
a_{j1} y_1 \partial_{y_1}  + 
a_{j2} y_2 \partial_{y_2} + \cdots + a_{jn} y_n \partial_{y_n} 
+ \gamma_j 
\quad \hbox{for} \quad
j = 1,2,\ldots,d .\nonumber
\end{eqnarray}

\smallskip

\noindent Here  $\,\partial_x^u = 
\partial_{x_1}^{u_1}\cdots \partial_{x_n}^{u_n}\,$
for $u \in \N^n$. In the notation of \cite{gkz89,gkz90,sst2},  
this is the system $H_A(-\beta,-\gamma)$ associated with 
the $(n+d) \times 2n$-matrix
\begin{equation}
\label{lawrencem}
A\quad :=\quad
\left(
\begin{array}{cc}
I_n & I_n\\
 0 \, 0 \, \cdots \, 0
& a_1 \, a_2 \, \cdots \, a_n \\
\end{array}
\right),
\end{equation}
where $I_n$ denotes the $n \times n$ identity matrix.
The matrix $A$ is called the {\it Lawrence lifting}
of $a_1,a_2,\ldots,a_n$. Such matrices
 play an important role in combinatorics
\cite[\S 9.3]{matroids} and Gr\"obner bases
 \cite[\S 7, page 55]{stu}.

We next introduce a combinatorial invariant associated
with a configuration of  vectors.
For the Lawrence lifting $A$, this invariant agrees with that of
the submatrix $ M := (a_1,\dots,a_n)$.  
The {\it matroid complex} of $M$ is the simplicial complex 
$\Delta(M)$
consisting of all subsets $I \subseteq \{1,\dots,n\}$ such 
that the
corresponding vectors $a_i, i\in I$ are linearly independent.
Let $\chi(M)$ denote the Euler characteristic of the 
matroid complex 
$\Delta(M)$, i.e. the sum of $(-1)^{|I|}$ for $I \in \Delta(M)$.
The integer $\chi(A) = \chi(M)$ equals
the M\"obius invariant of the dual matroid 
\cite[Proposition 7.4.7]{bjorner} and, via Zaslavsky's Theorem
\cite{zaslavsky}, it counts the regions of the
 hyperplane arrangement (\ref{CentralArr}).
Lemma \ref{nbc} implies
\begin{equation}
\label{atmost}
 |\chi(A)| \quad \leq \quad
 \binom{n-1}{d} , 
\end{equation}
with equality if
all $d$-tuples   $\{a_{i_1}, \ldots,a_{i_d}\}$
are linearly independent.

We note that $\,\chi(A) = 0 \,$
if and only if $A$ has a {\it coloop}, i.e., some linear
functional on $\R^d$ vanishes on all but one
of the points $a_1,\ldots,a_n$. If this is the case, 
then every $A$-hypergeometric function is a 
monomial times a solution of a smaller system
 (\ref{system}) gotten by {\it contracting
the coloops}. Thus, we will
assume without loss of generality that $\chi(A) \not= 0.$

A rational function $f$ in $x_1,\dots \! ,x_n$, 
$ y_1,\dots\!,y_n $
is called {\it unstable} if it is annihilated by some
iterated derivative $\,\partial_x^u\,\partial_y^v$.
Otherwise we say that $f$ is {\it stable}.  
Thus $f$
is unstable if it is a linear combination of rational functions
that depend  polynomially 
on at least one of the variables. 
We denote by 
$R(\beta,\gamma)$ the vector space of   rational
solutions of $H_{A}(-\beta,-\gamma)$, by
$U(\beta,\gamma)$ the subspace of unstable rational
solutions, and we set
\begin{equation}\label{stablequotient}
\ScS(\beta,\gamma)\ := \ R(\beta,\gamma)/U(\beta,\gamma)\,.
\end{equation}
Our main result gives an
integral representation for stable rational
$A$-hypergeometric functions, when $A$ is 
the Lawrence configuration (\ref{lawrencem}).

\begin{theorem}
\label{maintheorem}  Let  $\beta \in \Z_{>0}^n$ and
$\gamma \in \Z^d$. The space $\ScS(\beta,\gamma)$ 
of  stable rational
$A$-hypergeometric functions of degree $(-\beta,-\gamma)$ 
has dimension
$\, | \chi(A)| \,$
and is spanned by binomial residues $R_\Gamma(x ,y)$.
\end{theorem}

We illustrate this theorem with three examples. First consider
$d=1, n= 3, a_1 = a_2 = a_3 = 1 \in \Z^1$, 
$\beta_1 \! = \!  \beta_2 \! = \! \beta_3 = 1$, and $\gamma = 3$.
The Euler  characteristic is $\,\chi(A) = - 2$. 
The binomial residues are the integrals
\begin{equation}
\int_\Gamma 
 \frac{t^3}{(x_1+t y_1) (x_2+t y_2)(x_3+t y_3)} 
\frac{dt}{t}.
\end{equation}
 By integrating around the three poles
$t = - x_i/y_i$,  we obtain
$$ R_1  \quad = \quad  \frac{ x_1^2}
{(x_1y_2 -x_2y_1)(x_1y_3-x_3y_1)y_1} $$
$$ R_2  \quad = \quad \frac{ x_2^2}
{(x_3y_2-x_2y_3)(x_1 y_2-x_2y_1)y_2} $$
$$ R_3  \quad = \quad \frac{x_3^2}
{(x_2y_3-x_3y_2)(x_1y_3-x_3y_1)y_3} $$
These  residues form a solution basis for the
hypergeometric system
\begin{eqnarray*}
 H_A(-\beta, -\gamma) \, = &
 \{ \,
\partial_{x_1} \partial_{y_2}  -  
\partial_{x_2} \partial_{y_1}  ,\,
\partial_{x_1} \partial_{y_3}  -  
\partial_{x_3} \partial_{y_1}   ,\,
\partial_{x_2} \partial_{y_3}  -  
\partial_{x_3} \partial_{y_2} , 
\qquad \\
&  \!\!\!\!\!\!\!\!\!\!\!\!\!\!\!\!\!\!\!
 x_1 \partial_{x_1} +  y_1 \partial_{y_1} +  1, \,\,
x_2 \partial_{x_2} +  y_2 \partial_{y_2} +  1, \,\,
x_3 \partial_{x_3} +  y_3 \partial_{y_3} +  1,  
 \\ & 
  y_1 \partial_{y_1} +  
  y_2 \partial_{y_2} +  
  y_3 \partial_{y_3} +  3 \, \} .
\end{eqnarray*}
This is the {\it Aomoto-Gel'fand system}
for a $2 \times 3$-matrix, which is
holonomic of rank $3$; see \cite[\S 1.5]{sst2}. The space 
$\ScS(\beta,\gamma)$ of rational
solutions modulo unstable rational solutions
has dimension $\, 2 = | \chi(A)|\, $, since
$$ R_1 + R_2 + R_3 \,\, = \,\, \frac{1} {y_1 y_2 y_3} $$
contains no $x_i$ and is hence  unstable.
This identity expresses the fact that
the sum of  all local residues of a rational $1$-form 
over  $\CP^1$ is zero.

\smallskip

Our second example is the Lawrence lifting of the 
twisted cubic curve:
\begin{equation}
\label{sixbyeight}
d=2,\, n=4, \qquad \,
A \quad = \quad
  \left(
\begin{array}{cccccccc}
1 & 0 & 0 & 0 & 1 & 0 & 0 & 0 \\
0 & 1 & 0 & 0 & 0 & 1 & 0 & 0 \\
0 & 0 & 1 & 0 & 0 & 0 & 1 & 0 \\
0 & 0 & 0 & 1 & 0 & 0 & 0 & 1 \\
0 & 0 & 0 & 0 & 1 & 1 & 1 & 1 \\
0 & 0 & 0 & 0 & 0 & 1 & 2 & 3 \\
\end{array}
\right)
\end{equation}
Fix $\beta = (1,1,1,1)$ and $\gamma = (1,1)$.
The space of $A$-hypergeometric functions
is $10$-dimensional, and the subspace of rational solutions
is $3$-dimensional. A basis
for $R(\beta,\gamma)$ consists of the three binomial residues
$$
R_{23} \quad = \quad  \frac{x_2 y_3^2 y_2}{
(x_3^2 y_2 y_4 - x_2 x_4 y_3^2)(x_1  x_3 y_2^2 - x_2^2 y_1  y_3)}
$$
$$
R_{24} \quad = \quad \frac{y_4 y_2 
( x_2^2 y_3 y_1+x_1 x_3  y_2^2)}
{(x_3^2 y_2 y_4 - x_2 x_4 y_3^2)
( x_2^3 y_1^2 y_4 - x_1^2 x_4 y_2^3 )}
$$
$$
R_{34} \quad = \quad
\frac{ x_3 x_4 y_3^3  y_4}
{(y_2 y_4 x_3^2-y_3^2 x_2 x_4)(y_3^3 x_4^2 x_1-x_3^3 y_1 y_4^2)}
$$

\noindent Other residues can be computed by the
{\it Orlik-Solomon relations} (cf.~\S 5):
$$ 
R_{14} \,\, = \,\, -R_{24}-R_{34}, \quad
R_{13} \,\,= \,\,-R_{23}+R_{34}, \quad
R_{12} \,\,= \,\, R_{23}+R_{24}.
$$

\smallskip

For our third example take $\{a_1,\ldots, a_n \}$ to be the
positive roots in  the root system of type $A_{d}$.
This means $n= \binom{d}{2}$ and
(\ref{localresidue}) looks like
$$
 \int_\Gamma 
 \frac{
t_1^{\gamma_1}  \cdots t_{d-1}^{\gamma_{d-1}} }
{\prod_{1 \leq i < j \leq d }  (x_{ij} + 
t_i t_j^{-1} y_{ij})^{\beta_{1j}}}
\frac{dt_1}{t_1} \wedge \dots \wedge
\frac{dt_{d-1}}{t_{d-1}},
$$
where $t_d = 1.$
This is the {\it Selberg type integral} studied by
Kaneko \cite{kaneko} and many others; 
see  \cite[Example 5.4.7]{sst2}.
The holonomic rank of the associated $A$-hypergeometric system
equals $\,d^{d-2} $, the number of labeled trees on $d$ vertices.
The following explicit formula for the number of stable
rational hypergeometric functions of Selberg type is given
in \cite{nps}:
$$ |\chi(A_d)| \quad = \quad 
(d-2) \cdot \! \sum_{k = 0}^{[(d-3)/2]} \! \binom{d-3}{ 2k}
 (d-1)^{d-3-2k} \cdot \prod_{i=1}^k (2i-1).
$$

\smallskip 

This paper is organized as follows.
In \S 2 we examine  hypergeometric Laurent series solutions,
 and we derive the upper bound in Theorem
\ref{maintheorem}.
In \S 3 we establish the connection to toric geometry, by
expressing binomial residues as {\it toric residues}
in the sense of  Cox \cite{cox}; see also \cite{ccd, global}.  
Formulas and algorithms for computing binomial residues
are presented in \S 4. In \S 5, we complete the proof of
Theorem \ref{maintheorem}, and we prove
Conjecture 5.7 from our previous paper \cite{rhf}
in the Lawrence case.

\section{Laurent series expansions and Gale duality}

In this section we establish the upper 
bound in Theorem \ref{maintheorem}
for arbitrary rational $A$-hypergeometric functions. The Lawrence
hypothesis is not needed for this. The main idea is to look
at series expansions, which leads to counting cells in a
hyperplane arrangement. We fix
an arbitrary integer $ r \times s$-matrix $A$ of rank $r$
and  an integer vector $\alpha \in \Z^r$.

\begin{definition}\label{def:hypergeom}
\cite{gkz89,gkz90,sst2}.
The {\it $A$-hypergeometric system}
 is  the left ideal 
$H_A(\alpha)$ in the Weyl algebra
$ \C\langle x_1,\dots,x_s,\pd 1,\dots,\pd s\rangle$
 generated by the 
{\sl toric operators}
$\ \partial^u - \partial^v\ $, for 
$ u,v\in \N^s$ such that   $A\cdot u=A\cdot v$, 
and the {\sl Euler operators}
$ \,\sum_{j=1}^s a_{ij} x_j \partial_j - \alpha_i \,$ 
for $\, i=1,\dots,r $.
A function $f(x_1,\dots,x_s)$, holomorphic on an open
set $U\subset \C^s$, is said to be {\it $A$-hypergeometric of
degree $\alpha$} if it is annihilated by 
the left ideal $H_A(\alpha)$.
\end{definition}


A rational $A$-hypergeometric function admits 
Laurent series expansions convergent in a suitable open set.  
In the terminology of \cite{sst2} these are {\sl logarithm-free
hypergeometric series} with integral exponents.  We review their
construction and refer  to \cite[\S 3.4]{sst2}
for  proofs and details.  

Given a vector $v\in \C^s$, we define its 
{\sl negative support} by
$$
{\rm nsupp}(v)\ :=
\ \bigl\{\,i\in \{ 1,\dots,s\} \,:\, v_i\in\Z_{<0} \,\bigr\}
$$ 
A vector $v\in\C^s$ is said to have 
{\sl minimal negative support} if
there is no integer vector $u$ in the kernel of $A$
such that ${\rm nsupp}(u+v)$ is
properly contained in ${\rm nsupp}(v)$.
The following set of integer vectors,
$$
N_v\ :=\ \{u\in \ker_\Z(A) : {\rm nsupp}(u+v) = 
{\rm nsupp}(v) \} ,
$$
is used to define the formal Laurent series
\begin{equation}\label{hypseries}
\phi_v(x) \quad := \quad
\sum_{u\in N_v} \frac {[v]_{u_-}}{[v+u]_{u_+}} \cdot x^{v+u}
\end{equation}
$$
\hbox{where} \quad
[v]_{u_-} = \prod_{i:u_i<0} \prod_{j=1}^{-u_i} (v_i-j+1) \,\,\,
\hbox{and} \,\,\, [v+u]_{u_+} = 
\prod_{i:u_i>0} \prod_{j=1}^{u_i} (v_i+j) .$$
The Weyl algebra acts on formal 
Laurent series by multiplication and 
differentiation. The following is 
Proposition~3.4.13 in \cite{sst2}:

\begin{prop}\label{3.4.13} Let $\alpha = A\cdot v$.
The series $\phi_v(x)$ is annihilated by $H_A(\alpha)$
if and only if the vector $v\in\C^s$ has minimal negative support.
\end{prop}

In order to ensure that the $A$-hypergeometric series $\phi_v(x)$
have a common domain of convergence, we fix a generic weight
vector $\,w\in \R^s $. A vector $v\in \C^s$ is called an
{\it exponent for $H_A(\alpha)$ with respect to $w$} if
$v$ has minimal negative support and
\begin{equation}\label{exponent}
 A\cdot v = \alpha \qquad \hbox{and} \qquad
 \langle w,v \rangle  \ = \ \min
\{ \, \langle w,u \rangle  \,: \, u\in v + N_v \, \}.
\end{equation}
The following is a restatement of
\cite[Theorem 3.4.14, Corollary 3.4.15]{sst2}:

\begin{theorem}\label{3.4.14}
The set  $\, \bigl\{ \, \phi_v   \, : \, v \in \Z^s
\,\,\, \hbox{and $v$ is an exponent} \, \bigr\} \,$
is a basis for the space of hypergeometric 
functions of degree $\alpha$ admitting a 
Laurent expansion convergent in
a certain open subset $\, U_w \, $ of $\C^s$.
\end{theorem}

For a more precise description of hypergeometric 
Laurent series, we next introduce the  oriented hyperplane
arrangement defined by the Gale dual (or matroid dual) to $A$.
Set $m := s -r $ and   let $B$ be an integral
$s \times m$ matrix whose columns are a $\Z$-basis of
$\ker_\Z(A)$. The matrix $B$ has rank $m$ and $A\cdot B = 0$. 
Note that $B$ is
well-defined modulo  right multiplication by elements of 
$GL(m,\Z)$.  We identify $B$ with its set of row vectors,
and we call this configuration the {\it Gale dual} of $A$:
$$
B \ =\ \{b_1,\ldots,b_{s}\}\ \subset
\  \Z^m\,.
$$

Our assumption $\chi(A) \not= 0$ translates into
the condition $b_j \not= 0$ for all $j=1,\dots,s$.  
As remarked in the Introduction, the study of $A$-hypergeometric
functions, for arbitrary $A$,  easily reduces to this case.

Fix an exponent $v \in \Z^s$.
We identify the lattice $\Z^m$ with the sublattice
$\,image_\Z (B) + v = \ker_\Z (A) + v\,$ of $\Z^s$
via the affine isomorphism  
$\, \lambda \mapsto B \cdot \lambda + v $.
Under this identification,  the affine hyperplane
\begin{equation}\label{hyperpl}
 \{\lambda \in \R^m\ :\ \langle b_j,\lambda \rangle 
= -v_j\}\ \ \ \ 
\end{equation}
corresponds to the coordinate hyperplane $x_j = 0$ 
in $\,\ker_\Z (A) + v  \, \subset \, \Z^s $.
Let ${\HH}$ denote the arrangement in $\R^m $
consisting of the hyperplanes (\ref{hyperpl})
for  $\,j=1,\ldots,s$. We define the {\it negative support}
of a vector $\lambda $ in $\R^m$ as the negative support
of its image under the above isomorphism:
$$ {\rm nsupp} (\lambda) \quad := \quad
\bigl\{ \, j \in \{1,2,\ldots,s\}\,\,: \,\,
 \left\langle b_j, \lambda 
\right\rangle \ < \  -v_j\,\bigr\}.  $$
The set of points with the same 
negative support will be called a {\it cell}
of the hyperplane arrangement ${\HH}$. 
Note that our definition of cell
differs slightly from the familiar 
subdivision into relatively open
polyhedra by the hyperplanes in ${\HH}$. 
Our cells are unions of these:
they are also polyhedra but they are usually not relatively open. 

Consider the following attributes of a cell $\Sigma$ in ${\HH}$.
We say that:
\begin{itemize}
\item $\Sigma$ is {\it bounded} if $\Sigma$ 
is a bounded subset of $\R^m$.
\item $\Sigma$ is {\it minimal}  if
the set $\,\Sigma\cap\Z^s \,$ is nonempty and the
support of the elements in this set is minimal
with respect to inclusion.
\item $\Sigma$ is {\it $w$-positive}, for a given vector 
 $w$ on $\R^n$, if there exists a real number $\rho$ such that
$\,\langle w, \lambda \rangle \geq \rho \,$ 
for all $\lambda \in \Sigma $.
\end{itemize}

We can now rewrite the 
hypergeometric series (\ref{hypseries}) as follows:
\begin{equation}\label{sumcell}
\phi_\Sigma  \quad := \quad   \phi_v  \quad = \quad
\sum_{\lambda \in \Sigma \cap \Z^m}
 \frac {[v]_{(B\lambda)_-}}{[v+B\lambda]_{(B \lambda)_+}}
 \cdot x^{B \lambda + v} 
\end{equation}
If $\Sigma$ is bounded then $\phi_\Sigma$ 
is a Laurent polynomial,
and if $\Sigma $ is $w$-positive then $\phi_\Sigma$ lies in
the Nilsson ring (cf.~\cite[\S 3.4]{sst2})  associated with $w$,
and hence defines an $A$-hypergeometric function on $U_w$
when $\Sigma$ is minimal.

The following is an immediate consequence of 
Theorem~\ref{3.4.14}:

\begin{prop}\label{laurentpolyn} 
The series $\phi_\Sigma$ where $\Sigma$ runs over all
 $w$-positive minimal cells in ${\HH}$ form
a basis for the space of $A$-hypergeometric 
functions of degree $\alpha$ admitting a Laurent 
expansion convergent in
 $\, U_w \, \subset \, \C^s$.
Restricting to bounded cells $\Sigma$, we get a basis
for the subspace of hypergeometric Laurent polynomials.
\end{prop}

Recall that a rational function $f$ in 
is called {\it unstable} if  there exists $u \in \N^s$ such that
the partial derivative $\partial^u(f)$ is identically zero. 

\begin{lemma}
\label{poscoor}
$A$-hypergeometric Laurent polynomials are unstable.
\end{lemma}

\proof
By  Proposition \ref{laurentpolyn}, it is enough to show
that  for any  bounded minimal chamber $\Sigma,$ the common
negative support of all monomials in $\phi_\Sigma$ does not
equal  $ \{1,\dots,s\}.$  In fact, suppose
$$
\Sigma \quad = \quad
\bigl\{ \lambda \in \R^m \,: \,
 \left\langle b_j, \lambda \right\rangle \ < \  -v_j\,\,\,
\hbox{for all} \,\, j =1,2,\ldots,m \,
\bigr\}.  $$
The negative support of any  lattice point  in 
$\Z^m \backslash \Sigma$
is a proper subset
of $\{1,2,\ldots,n\}$, and we conclude that 
$\Sigma$ is not minimal.
\ \qed

\smallskip

If we differentiate an $A$-hypergeometric 
function of degree $\alpha$
with respect to $x_i$ then we get an $A$-hypergeometric function
of degree $\alpha - a_i$. If we iterate this process long enough,
for all variables, then only the stable functions survive. The 
following definition is intended to make this more precise.
The {\it Euler-Jacobi cone} is the  open cone in $\R^s$:
$$-{\rm Int}({\rm pos}(A)) 
\  = \  \bigl\{
\nu_1 a_1 + 
\nu_2 a_2 +  \cdots +  \nu_s a_s \,  :   \,\,
\nu_i \in \R_{< 0} \  \hbox{ for all } \, i \, \bigr\}. $$
Note that $(-\beta,-\gamma)$ lies 
in the Euler-Jacobi cone in Example (\ref{sixbyeight}).

\begin{proposition}
\label{ejprop}
If $\,\alpha \in -{\rm Int}({\rm pos}(A)) $, then
every $A$-hypergeometric series of degree $\alpha$ is stable.
\end{proposition}

\proof 
If suffices to show that none of the hypergeometric series
$\phi_v$ is unstable.
Fix a strictly negative vector $\nu \in \Q_{<0}^s$
with $\, A \nu = A v = \alpha$. Let $k$ be a positive integer
such that
$k\nu \in\Z^s$. For each integer $\ell \in \N$, the vector
$\,v + \ell(v  - \nu) \,$ has negative support 
contained in $\nsupp(v)$.
Since $v$ is minimal, we conclude that
$\nsupp (v + \ell(v  - \nu)) = \nsupp(v)$
for all $\ell \in N$. 
Let $I := \{ i \in \{1,\dots,s \} \, : \, v_i \geq 0 \}.$ 
For all $i \in I,$ we have $v_i > \nu_i$,  and so
all the coordinates in $I$ of the vectors 
$\,v + \ell(v  - \nu)\,$
strictly increase with $\ell$. This shows that
$\phi_v$ cannot be decomposed as a finite sum of 
Laurent series that
depend polynomially on one variable. 
\ \qed

\smallskip

\begin{theorem}
\label{bound}
If $\,\alpha \in -{\rm Int}({\rm pos}(A)) $, then 
the dimension of the  space of $A$-hypergeometric Laurent series
of degree $\alpha$ with a common domain of 
convergence is bounded above
by the Euler characteristic $| \chi(A) |$.
\end{theorem}

\proof
Consider the central hyperplane arrangement 
gotten from ${\HH}$
by translating all $s$ hyperplanes so as to 
pass through the origin.
This central arrangement consists of the $s$ hyperplanes
\begin{equation}
\label{CentralArr}
 \bigl\{ \, \lambda \in \R^m \,\,: \,\,
\langle b_j, \lambda \rangle \, = \, 0 \,\bigr\} 
\qquad \hbox{for} \,\,\, j = 1,2,\ldots,s . 
\end{equation}
Since $\alpha$ is in the Euler-Jacobi cone,
the minimal cells $\Sigma$ of ${\HH}$ are all unbounded
and correspond to certain maximal cones of the 
central arrangement (\ref{CentralArr}). Fix a generic
linear functional $w$ on $\R^m$. 
A basis for the relevant space of
$A$-hypergeometric Laurent series is indexed by the
$w$-bounded, minimal cells of ${\HH}$. Their number is
bounded above by the  number of $w$-bounded maximal cones in the
central arrangement.

A classical result in  combinatorics due to Zaslavsky
\cite{zaslavsky} states
that the number of $w$-bounded maximal cones is
the absolute value of the  M\"obius invariant $\mu(B)$ of 
the matroid associated with $B$. Our assertion now follows
from the following identity from
\cite[Proposition 7.4.7 (i)]{bjorner}:
\begin{equation}
\label{MUisCHI}
 | \mu(B) | \quad = \quad | \chi(A) | .
\end{equation}
In words, the M\"obius invariant of a matroid equals
(up to sign) the Euler characteristic of the dual matroid.
{\qed}

\begin{corollary} 
\label{lowb}
For any $\alpha \in \Z^d$, the complex vector space of  rational
$A$-hypergeometric functions of degree $\alpha$ modulo the subspace of
unstable functions has dimension at most $|\chi(A)|$.
\end{corollary} 

\proof
We represent the rational $A$-hypergeometric functions
by Laurent series expansions which have a common domain
of convergence. Hence it suffices to prove the asserted dimension 
bound for the space of convergent $A$-hypergeometric Laurent series
modulo unstable ones.

Choose $u \in \N^s$ so that  $\alpha - A u $ lies in the Euler-Jacobi cone.
The  operator $\partial^u$ induces a monomorphism
from $\ScS(\alpha)$ into $\ScS(\alpha - A u)$.  By
Proposition~\ref{ejprop},  $\ScS(\alpha - A u) \cong R(\alpha -
A u)$, hence the dimension bound follows from
Theorem~\ref{bound} applied to
$\alpha - A u$.
{\qed}

\vskip .1cm

Passing from $\{a_1,\ldots,a_n\}$ to its Lawrence
lifting (\ref{lawrencem}) corresponds under Gale duality to 
the operation of replacing $ \{b_1,\ldots,b_n\}$ by its 
symmetrization $\{b_1,\ldots,b_n,-b_1,\ldots,-b_n \}$; 
see  \cite[Proposition 9.3.2]{matroids}. This process does not
change the geometry of the hyperplane arrangement
(\ref{CentralArr}) and hence
it does not change the M\"obius invariant $\mu(B)$. In view of 
(\ref{MUisCHI}), we conclude that the Euler characteristic of
$\{a_1,\ldots,a_n\}$ equals the Euler characteristic of its
Lawrence lifting  as stated in the Introduction. 
Corollary \ref{lowb} implies
the upper bound in Theorem \ref{maintheorem}.

\begin{corollary}
The space  $\ScS(\beta,\gamma)$
has dimension at most $\, | \chi(A)| \,$
\end{corollary}

\smallskip
We conclude this section with one more result
from matroid theory which we need to complete the proof
of Theorem~\ref{maintheorem}.  
A maximally independent subset of $B$ is a {\it basis} of $B$.
Note that 
$\,\{b_j : j \in J \}\,$ is a basis of $B$ if and only if
$\,\{a_j : j \not\in J \}\,$ is a basis of $A$.
A minimally-dependent subset of $B$ is 
a {\it circuit} of $B$. If 
$\,C = \{b_{i_1},\ldots,b_{i_t}\}\,$ is a circuit
and $i_1  < \cdots  < i_t$ then the set
$\,C \backslash \{b_{i_t} \}\,$ is a {\it broken
circuit}. A basis of $B$ 
is called an {\it nbc-basis} if it contains no
broken circuits.

\begin{lemma}\label{nbc}
The number of nbc-bases of $B$ equals $|\chi(A)|$.
\end{lemma}

\proof
This result follows from (\ref{MUisCHI}) and
Proposition 7.4.5 in \cite{bjorner}.
\qed

\section{Binomial residues and toric geometry} 

This section is concerned with global residues of meromorphic
forms whose polar divisor is a union of hypersurfaces defined
by binomials. The analogous case when the polar divisor is
defined by linear forms has been extensively studied,  for instance, 
by Varchenko \cite{var} and Brion-Vergne \cite{brion}. Our
situation can be regarded as a multiplicative analogue to that
theory. The binomial hypersurfaces are embedded in a suitable 
projective toric variety, which places binomial residues into the
framework of  toric residues  \cite{ccd, global, cox}.
This will allow us in \S 5 to find bases of $A$-hypergeometric stable
rational functions for Lawrence liftings in terms of binomial residues, and to
give a geometric meaning to the linear dependencies among
binomial residues. We refer to \cite{griffiths, tsikh} for the
definition and basic  properties of   Grothendieck
residues.

Let $X$ be a complete $d$-dimensional toric variety and $S$ its
homogeneous coordinate ring in the sense of Cox \cite{coxhomog}. 
Homogeneous polynomials in $S$ may be thought
of as sections of  coherent sheaves over $X$ and, 
consequently, their zero-loci are well defined divisors in $X$.
Let $T \simeq (\C^*)^d$ denote the dense torus in $X$.
Suppose  $G_0,G_1,\ldots,G_d$ are homogeneous polynomials in
$S$ whose divisors $D_i$ satisfy
\begin{equation}\label{emptyy}
D_0\cap D_1 \cap \cdots\cap D_d = \emptyset.
\end{equation}
Any homogeneous polynomial $H$ of critical degree determines
a meromorphic $d$-form on $X$ with polar divisor
 contained in $D_0\cup \cdots\cup D_d$,
$$\Phi(H) \quad = \quad \frac  {H\,\Omega_X}{ G_0 G_1 \cdots G_d}\, ,
$$
where $\Omega_X$ is a choice of an Euler form on $X$ \cite{batyrevcox}.  
The $d$-form $\Phi(H)$ defines a \v Cech cohomology
class $[\Phi(H)] \in H^d (X,\widehat\Omega_X^d)$ relative to the open cover 
$\{ X \backslash D_i\}_{i=0,\dots,d }$ of $X$.  Here
$\widehat\Omega_X^d$ denotes the sheaf of Zariski $d$-forms on $X$.
The class $[\Phi(H)]$ is alternating with respect to permutations
of $G_0,\dots, G_d$. If $H$ lies in the ideal 
$\langle G_0,\ldots,G_d\rangle$ of $S$ then $\Phi(H)$ is a \v Cech coboundary.
Thus, $[\Phi(H)]$ depends only on the image of
the polynomial $H$ in  the quotient ring
$\,S/\langle G_0,\ldots,G_d \rangle$.  

The {\it toric residue\/}
$\rx_G(\Phi(H)) \in \C$ is given by the formula
$$
\rx_G(\Phi(H))  \,\,= \,\, {\rm Tr}_X([\Phi(H)]),
$$
where $\,{\rm Tr}_X : H^d(X,\widehat\Omega_X^d) \to \C \,$ is the trace map.

The following proposition can be deduced from Stokes Theorem 
(cf. \cite{griffiths}, \cite[\S 7.2]{tsikh}). 
It follows directly from the definition of  toric residue.

\begin{proposition}
\label{stokes}
If the  polar locus of
the $d$-form $\Phi(H)$ is contained in the union of
only $d$ divisors, say
$D_1\cup \cdots\cup D_d$, then $\rx_G(\Phi(H)) = 0$.
\end{proposition}

The relationship between toric residues and the usual notion of
multidimensional residues is given by the following result.

\begin{theorem}\label{localglobalth}
Let $G_0,\dots,G_d\in S$ satisfy (\ref{emptyy}) and suppose
\begin{equation}\label{finite}
V^{0} := D_1\cap \cdots\cap D_d \,\, \subset \,\, T
\end{equation}
Then 
\begin{equation}\label{localglobaleq}
\rx_G(\Phi) \ =\ \sum_{\xi \in V^{0}} {\rm Res}_\xi(\Phi|_T)
\end{equation}
where ${\rm Res}_\xi(\Phi|_T)$ denotes the (local) Grothendieck
residue at $\xi$ of the meromorphic form $\Phi$
restricted to the torus and relative to the divisors 
$D_1\cap T, \dots,D_d\cap T$. 
\end{theorem}

\proof
We note, first of all, that (\ref{emptyy}) implies that 
$V^{0}$ is a finite set and hence the sum in 
 (\ref{localglobaleq}) makes sense.  Moreover, as shown
in  \cite[\S II.7.2]{tsikh}, the local residues in the
right-hand side of (\ref{localglobaleq})  depend only on the
divisors $D_i\cap T$ and not on the choice of local defining
equations.

If $X$ is simplicial, 
then (\ref{localglobaleq})  is the content of Theorem 0.4
in \cite{ccd}. For general $X$ we argue as in the proof of Theorem
4 in  \cite{global}.  \ \qed 

\smallskip

We consider now the binomial case which is relevant in this paper.
Let $a_1, \dots,a_n \in \Z^d$ as in the Introduction.
Let $\Delta_i$ denote the segment
$[0,a_i] \subset \R^d$ and  $\Delta = \Delta_1+\cdots+\Delta_n$
their Minkowski sum. This is a {\it zonotope}, that is, a polytope
all of whose faces are centrally symmetric \cite[\S 2.2]{matroids}.
Let $\eta_1,\ldots,\eta_{2p}$ denote the inner normals of the
facets of the zonotope $\Delta$, where  $\eta_j = -\eta_{p+j}$.
 We can write
$$\Delta \ =\ \bigl\{\,m\in\R^d \,\, : \,\, \langle m, \eta_j\rangle
\, \geq \! \sum_{i: \langle \eta_j, a_i\rangle < 0 } \!
  \langle \eta_j,  a_i\rangle 
\,; \, j=1,\dots,2p \,\bigr\}$$
We consider the associated projective toric
variety $X_\Delta$.   
The homogeneous coordinate ring of $X_\Delta$ is 
the polynomial ring $S = \C[z_1,\dots,z_{2p}]$.
The
monomials $\,t_j \,:=\, \prod_{i=1}^{2p}
\bigl(z_i^{\eta_{ij}}\bigr)\,$, for $j = 1,2,\ldots,d$,
have degree zero and define coordinates in the torus
$T\subset X_\Delta$.

To each binomial $f_i := x_i + y_i t^{a_i}$ 
in the denominator of the kernel of (\ref{localresidue})
we associate the homogeneous polynomial
$$F_i (z) \ :=\ 
 x_i  \prod_{\langle \eta_j, a_i\rangle < 0} z_j^{- \langle
\eta_j, a_i\rangle }+ y_i 
\prod_{\langle \eta_j, a_i\rangle > 0} z_j^{\langle \eta_j, 
a_i\rangle }.$$
The divisor $Y_i := \{F_i(z) =0\}\subset X_\Delta$ is the closure
of the divisor $\{f_i(t) =0\}\subset T$.  Moreover, for
$\beta\in \Z_{>0}^n$ and $ \gamma\in \Z^d$,
the  $d$-form on $T$,
\begin{equation}\label{meromtorus}
 \phi (\beta,\gamma) \quad = \quad
\frac{t^{\gamma}} 
{ f_1^{\beta_1} 
\cdots f_n^{\beta_n}} \frac{dt_1}{t_1} \wedge \dots \wedge
\frac{dt_d}{t_d}, 
\end{equation}
extends to the following meromorphic $d$-form on the toric variety $X_\Delta$:
\begin{equation}\label{meromtoric}
\Phi (\beta,\gamma) \quad  = \quad \frac{z^{h(\beta,\gamma)}} 
{ F_1^{\beta_1} 
\cdots F_n^{\beta_n}} \, {\Omega_\Delta},
\end{equation}
$$ 
\!\! \hbox{where} \quad
h_j(\beta,\gamma) \,\,\, = \,\,\, \langle \eta_j, \gamma \rangle - 
\!\! \sum_{\langle
\eta_j, a_i\rangle < 0}
  \langle \eta_j, \beta_i a_i\rangle -1,\ \  j=1,\dots,2p.
$$
The polar divisor of $\Phi (\beta,\gamma)$ is the union
of the divisors $Y_1,\ldots,Y_n$ and coordinate divisors
$\,\{z_\ell = 0\}$ for indices $\ell$ with
$h_\ell(\beta,\gamma) < 0$. 
For degrees in the Euler-Jacobi cone
such indices $\ell$ do not exist.  Indeed, 
 \begin{equation}
\label{eje} -{\rm Int}({\rm pos}(A))=
\bigl\{  (-\beta, -\gamma)  \in \R^{n+d}
 \,: \,  \beta_i > 0
; \, 
  h_j(\beta,\gamma) +1 > 0 \bigr\} 
\end{equation}
Thus, if $(-\beta,-\gamma)$ lies in the
Euler-Jacobi cone, 
the polar divisor of $\Phi (\beta,\gamma)$ equals
$Y_1 \cup \cdots \cup Y_n$.

\smallskip

We are now prepared to give a precise definition of binomial
residues. Fix an index set $I = \{ 1\leq i_1<\cdots < i_d \leq n\}$ 
such that the corresponding vectors 
$a_i$, $i\in I$, are linearly independent.  For $k=1,\dots,d$, set 
$G_k^I = F_{i_k}$ and $D_k = \{G_k^I = 0\}$. For generic values of the
coefficients 
$x_i, y_i$, $i\in I$, the divisors $D_1,\dots,D_d$ satisfy
(\ref{finite}). 

\begin{definition}
For $\beta\in \Z_{>0}^n$ and $ \gamma\in \Z^d$, 
 let
$$G_0^I \ =\ \bigl(\prod_{\ell:h_\ell(\beta,\gamma) < 0} \!\!\! z_\ell
\, \bigr)\ 
\cdot\ \bigl( \prod_{j\not\in I} F_j\bigr). $$
Define the following quantity which depends on
$\,x_1,\ldots,x_n,y_1,\ldots,y_n $:
$$
R_I(\beta,\gamma) \quad := \quad
\rx_{G^I}(\Phi (\beta,\gamma)).
$$
\end{definition}

\smallskip

Each local residue in the right-hand side of (\ref{localglobaleq})
may be written as an integral over a $d$-cycle ``around"  the point
$\xi \in V^0$.  Since for generic values of the coefficients, the map
$f_I = (f_{i_1},\dots,f_{i_d}) \colon T \to \C^d$ is proper, 
it follows from \cite[\S II.8]{tsikh} that the total sum of
residues (\ref{localglobaleq})
may be written as a single integral,

\begin{equation}
\label{singleintegral}
R_I(\beta,\gamma) \quad = \quad  \left(\frac 1 { 2 \pi i} \right)^{\! d}
\int_{\Gamma(I,x,y)} \frac{t^{\gamma}} 
{ f_1^{\beta_1} 
\cdots f_n^{\beta_n}} \frac{dt_1}{t_1} \wedge \dots \wedge
\frac{dt_d}{t_d},
\end{equation}
where $\Gamma(I,x,y)$ is the compact real
$d$-cycle  $\Gamma(I,x,y) \subset T$ defined by
$\{ |f_{i_1}| = \varepsilon_1, \dots, |f_{i_d}| = \varepsilon_d
\} $ for small positive $\varepsilon_1, \dots, \varepsilon_d.$
Moreover,  the cycle  $\Gamma(I,x,y)$ 
can be locally replaced by a cohomologous cycle
$\Gamma(I)$ independent of $(x_1,\dots,x_n,y_1,\dots,
y_n)$. See  \cite[\S 5.4]{sst2} for further details.

We close this section with the observation that the
``basic binomial residue''  $\,R_I(\beta,\gamma)\,$  
is indeed a rational $A$-hypergeometric  function.

\begin{lemma}
\label{sososos}
The toric residue
$R_I(\beta,\gamma)$ is a rational function
of  $(x,y) $ and is annihilated 
by the hypergeometric system (\ref{system}).
\end{lemma}

\proof
For any choice of polynomials $G_0,\ldots,G_d$, the
trace map $\,{\rm Tr}_X\,$ in the definition of the toric
residue has its image in the subfield of $\C$
generated by the coefficients of the $G_i$. This implies
that $R_I(\beta,\gamma)$ is an element 
in the rational function field
$\Q( x_1,\ldots,x_n,y_1,\ldots,y_n )$.

The kernel of the integral (\ref{sososos}) is annihilated by
the toric operators $\,\partial_x^u\,\partial_y^v\, 
-\,\partial_x^v\,\partial_y^u \,$ in (\ref{system}).
Hence so is the integral itself,
by diffentiating under the integral sign.
Specifically, it  follows from \cite[Lemma 6]{global} that 
\begin{equation}\label{derx}
\partial_{x_i} R_I(\beta,\gamma)\  = \ 
 - \beta_i \, R_I(\beta + e_i,\gamma), \quad \hbox{and}
\end{equation} 
\begin{equation}\label{dery}
\partial_{y_i} R_I(\beta,\gamma)\  = \  - \beta_i \, 
R_I(\beta + e_i,\gamma+a_i),
\end{equation}  
where $e_1,\dots,e_d$ is the standard basis of
 $\R^d$. The verification of the homogeneity
equations is immediate from the expression (\ref{meromtorus}) for the form
$\phi(\beta,\gamma)$.  Hence
$R_I(\beta,\gamma)$ is a rational
solution of $H_A(-\beta,-\gamma)$.  
 \ \qed

\section{Computing binomial residues}

In this section we present methods for computing 
the binomial residue $R_I(\beta,\gamma)$. Here  $I = \{i_1,\ldots,i_d\}$ 
is a fixed column basis of the matrix $M = (a_1,\dots,a_n)$. Let
$M_I$ denote the non-singular $d\times d$ matrix with columns
$a_i$, $i\in I$. Write $ M_I^{-1} = (\mu_{ij}) \in GL(d,\Q)$.
We set $\,V_I \,= \, \{ \xi \in T : f_i( \xi )=0 \,$ for all $\, i\in I\}$.
The points in $V_I $ are in bijection with the characters 
$\theta\in {\rm Hom}(\Z^d, \C^*)$  satisfying
 $\theta(a_i) = -1$, for all $i\in I$.  
The point $\xi^\theta = (\xi^\theta_1,\dots,\xi^\theta_d) \in V_I$ 
indexed by $\theta$ has coordinates
$$\xi_j^\theta \ =\ \theta(e_j)\cdot\prod_{i\in I} \left(\frac
{x_i}{y_i}\right)^{\mu_{ij}}$$
There are  $\det(M_I)$-many simple roots $\xi^\theta$
provided all $x_i, y_i$ are nonzero.

Let $g$ be a function  meromorphic on the torus 
$T = (\C^*)^d$ and regular at a simple root $\xi \in V_I$. Then 
the {  local Grothendieck residue} of  the 
meromorphic $d$-form $ \,
 \frac{g} { f_{i_1} \cdots f_{i_d}} \frac{dt_1}{t_1} \wedge \dots \wedge
\frac{dt_d}{t_d}\,$ at the point $\xi$ equals
\begin{equation}\label{simpleresidue}
R_{I,\xi}[g] \quad = \quad 
\frac {g(\xi)}{J_I(\xi)}
\end{equation}
where $J_I$ denotes the {\it toric Jacobian} of
the binomials $f_i = x_i + y_i t^{a_i}$:
$$J_I(t)\ =\ \det\left( t_j \frac {\partial f_i}{\partial t_j}
\right)_{i\in I}^{ j = 1,\ldots,d }
\ =\ 
\det M_I\cdot (\prod_{i\in I} y_i) \cdot 
{t^{a_I}} . $$
Here $a_I =a_{i_1}+\cdots+a_{i_d}$. 
We deduce the following identity
\begin{equation}\label{jacobian}
J_I(\xi)\ =\ (-1)^d\cdot
\det M_I\cdot (\prod_{i\in I} x_i) \qquad
\hbox{for all} \,\,\,
\xi \in V_I . 
\end{equation}
We obtain the following procedure for summing
(\ref{simpleresidue}) over all $\xi \in V_I$.

\begin{algorithm}
\label{firstalgorithm}
 {\sl (Computing global residues
using Gr\"obner bases)}
\hfill

\noindent {\sl Input:} 
A $d \times d$-integer matrix $M_I$ of rank $d$,
a Laurent polynomial $\,g(t) $.

\noindent {\sl Output:} The {\sl global residue}
$$R_{I} [g]: = \sum_{\xi\in V_I} R_{I,\xi}[g] $$

\begin{itemize}
\item[(1)]
 Fix the field $K = \Q(x_1,\ldots,x_n,y_1,\ldots,y_n)$
and write the Laurent polynomial ring over $K$ as
a quotient of a polynomial ring:
$$ \quad \qquad K[t_1,\ldots,t_d, t_1^{-1},\ldots, t_d^{-1}] \quad  = \quad
K[t_0,t_1,\ldots,t_d]/\langle t_0 t_1 \cdots t_d - 1 \rangle .$$
\item[(2)]
Compute any Gr\"obner basis $G$ for its ideal
$\,\langle f_{i_1},\ldots,f_{i_d} \rangle $.
\item[(3)]
Let $B$ be the set of 
standard monomials for $G$ in
$K[t_0,\ldots,t_d]$
\item[(4)] Compute the {\it trace} of $g$ modulo $B$ as follows:
$$
\sum_{\xi \in V_I} g(\xi)
\quad = \quad
\sum_{t^b \in B}
{\rm coeff}_{t^b} \bigl( {\rm normalform}_G ( t^b \cdot g (t) ) \bigr)
$$
\item[(5)] Output the result of step (4) divided by the monomial 
in (\ref{jacobian}).
\end{itemize}
\end{algorithm}

\smallskip

The output produced by the above algorithm is
a rational function in $x_i,y_i$ and the coefficients of $g$.
In the case when $g$ is a Laurent monomial, one can give
a completely explicit formula for that output.

\begin{lemma}
\label{computation}
Let  $\gamma\in\Z^d$. If 
$\,\nu =  M_I^{-1} \cdot \gamma \,$ lies in the lattice $\Z^d $ then
\begin{equation}\label{glorestq}
R_{I}[t^\gamma](x,y) \ =\  \frac {(-1)^{|\nu|+d}}{\det(M_I)}
\cdot \prod_{i\in I} x_i^{\nu_i -1}\,y_i^{-\nu_i }\,.
\end{equation}
Otherwise the {\sl global residue}
$\,R_{I} [t^{\gamma}] \,$ is zero.
\end{lemma}

\proof 
It follows from (\ref{simpleresidue}) and (\ref{jacobian}) that
$$
R_{I,\xi^\theta}[t^\gamma](x,y) \ =\ \theta(\gamma) \cdot \frac
{(-1)^d}{\det(M_I)}
\cdot \prod_{i\in I} x_i^{\nu_i -1}\,y_i^{-\nu_i }$$
where $\nu_i :=\sum_{j=1}^d \mu_{ij}\gamma_j$, $i\in I$.
Thus, the global residue
 is given by
$$
R_{I}[t^\gamma](x,y) \ =\ \left(\sum_\theta \theta(\gamma)\right) \cdot
\frac {(-1)^d}{\det(M_I)}
\cdot \prod_{i\in I} x_i^{\nu_i -1}\,y_i^{-\nu_i }$$
and consequently it vanishes unless 
$\gamma\in M_I\cdot \Z^d$.  In this case we have
(\ref{glorestq}) for
$\ \gamma=\sum_{i\in I} \nu_ia_i\,$, and $\ |\nu| := \sum_{i\in I} \nu_i$.
\ \qed 

\smallskip

We  now   compute the binomial residue
$R_I(\beta,\gamma)$ for $I = \{ i_1, \ldots, i_d \}$ as above.
In view of (\ref{derx}) and (\ref{dery}),
it suffices to consider the case $\beta= \one := (1,\dots,1)$.
Set $\, J \, :=  \{1,\ldots,n \} \backslash I $ and
let $M_J$ denote the matrix whose columns are the vectors $a_j, j\in J$.
Since the coefficients are generic, none of the polynomials
$f_j$, $j\in J$ vanishes on any point of $V_I$ and hence 
\begin{equation}
\label{effjay} 
R_I(\one,\gamma) \, = \, R_I[\,t^\gamma/f_J(t)\,](x,y) \qquad \hbox{where} 
\quad f_J(t) = \prod_{j \in J} f_j(t) .
\end{equation}
This gives rise to the following symbolic algorithm for binomial residues.

\begin{algorithm} {\sl (Computing binomial residues)}
\hfill

\noindent {\sl Input:} Vectors $a_1,\ldots,a_n$
and $\gamma$ as above, and a basis $I = \{i_1,\ldots,i_d\}$.

\noindent {\sl Output:} The rational function
$\,R_I(\one,\gamma) \,$ of $x_1,\ldots,x_n,y_1,\ldots,y_n$.

\begin{itemize}
\item[(1)] Run steps (1), (2) and (3) of Algorithm \ref{firstalgorithm}.

\item[(2)] Using linear algebra over the field $K$, compute the unique
 polynomial $\,g(t) \, = \,\sum_{t^b \in B} c_b \cdot t^b \,$ 
such that all $c_b $ lie in $K$ and
$\, g(t) \cdot  f_J(t) - t^\gamma \,$ reduces to zero modulo 
the Gr\"obner basis $G$.
\item[(3)] Run steps (4) and (5) of Algorithm \ref{firstalgorithm}.
\end{itemize}
\end{algorithm}

\smallskip

The output of this algorithm is an element of the field $K$.
It is nonzero and has the following expansion as a  Laurent series
in $x_i, y_i$.

\begin{proposition}\label{globalri}
Suppose $\gamma\in M\cdot\Z^n$.
Then $\,R_I(\one,\gamma)\not= 0\,$ and
\begin{equation}\label{gloresexp}
R_I(\one,\gamma)\ =\ 
\frac {1}{\det(M_I)}\ \sum
 (-1)^{d + |\nu| + |\mu|}\cdot \prod_{i\in I}  \frac
{x_i^{\nu_i-1}}{y_i^{\nu_i}}\cdot 
\prod_{j\in J} \frac {y_j^{\mu_j}}{x_j^{\mu_j+1}},
\end{equation}
where the sum is over $\nu\in\Z^I$ and $\mu\in \N^{J}$ such that
$\, M_I\cdot\nu - M_J\cdot \mu = \gamma $.
Moreover, for every $\beta\in \Z^n_{>0}$, the residue
$R_I(\beta,\gamma)$ is a stable rational hypergeometric
function.
\end{proposition}

\proof 
  We expand
\begin{equation}\label{expansion}
 t^\gamma\cdot f_J(t)^{-1} \ =\  \sum_{\mu\in \N^J} \  \prod_{j\in J}
\left( {y_j^{\mu_j}}\cdot{x_j^{-\mu_j-1}} \right)\ t^{\gamma
+M_J\cdot
\mu}\,.
\end{equation}
Applying (\ref{glorestq}) to each term of
(\ref{expansion}) yields the  Laurent expansion (\ref{gloresexp}).

Suppose now that $\gamma = \La_I\cdot\nu_0 - \La_J\cdot\mu_0$,
$\nu_0 \in \Z^I$,  $\mu_0 \in \Z^J$.
There exists a vector  
${m}\in \Z_{>0}^J$ such that 
$m_j a_j \in  M_I\cdot \Z^d$.  Hence for $k\in\N$,
$$\gamma + \La_J\cdot(\mu_0 + k m) \in M_I\cdot \Z^d$$
and $\mu_0 + k m$ is non-negative for $k \gg 0$.
 Hence, the series 
(\ref{gloresexp}) contains infinitely many non-zero terms.
This shows that $R_I(\one,\gamma) \not= 0 $.

Suppose now that $\beta\in \Z_{>0}^n$ is arbitrary.  In view of
(\ref{derx}), it suffices to show that the derivative
$\partial_x^{\beta - \one}$ of the series (\ref{gloresexp})
contains infinitely many powers of each of the variables 
$x_\ell$, $\ell=1,\dots,n$.  The previous argument shows that this is
indeed the case for $x_j$, $j\in J$ and also for a variable
$x_{i_0}$, $i_0\in I$, unless every vector $a_j$, $j\in J$, is in the
$\Q$-span of $\{a_i, i\in I, i\not= i_0\}$.  But this would mean that
the points $a_k$, $k\not= i_0$ would define a coloop in $A$ which
is impossible by assumption.
\ \qed 

\smallskip

Our final task in this section is to identify the irreducible
factors in the denominators of these binomial residues.
Let $C \subseteq \{1,\ldots,n\}$ be a {\it circuit}, i.e.,
the set $ \{a_i, i\in C\}$ obeys a unique (up to sign)
linear relation  $\sum_{i\in C} m_i a_i =0$
over $\Z$ such that ${\rm gcd}(m_i , \, i \in C) =1$. Then
$${\rm Res}(C;x,y) \quad = \quad
\prod_{m_i>0} x_i^{m_i} \prod_{m_j <0}
y_j^{m_j}
- (-1)^{|C|} \prod_{m_i>0} y_i^{m_i} \prod_{m_j <0} x_j^{m_j}$$
is the {\sl resultant} of the binomials $f_i, i\in C$.
In fact, the singular locus of $H_A(-\beta,-\gamma)$ is described
by the product of all the variables and all the resultants ${\rm Res}
(C;x,y)$ as $C$ ranges over the circuits (cf. \cite{rhf},\cite{gkz89}). 
Let $I$ be a basis as above. Note that for each $j\not\in I$,
there exists a unique subset $I'(j)\subseteq I$, such that 
$I(j) := I'(j)\cup\{j\}$ is a circuit.  

\begin{theorem}\label{formula2} The binomial residue,
defined by $I, \beta , \gamma$ as above, equals
\begin{equation}
R_I(\beta,\gamma)\ =\ 
\frac {P(x,y)} {x^a\,y^b\,\prod_{j\not\in I} {\rm Res}(I(j);x,y)^{c_j}}
\qquad \hbox{with all $c_j > 0$}
\end{equation}
where $P(x,y)$ is a  polynomial relatively prime from the 
denominator.
\end{theorem}

 \proof  We may assume that $\beta =  \one$.
It follows
from a variant of
Theorem~1.4 in \cite{tokyo} that $R_I(\one,\gamma)$ is a rational
function whose denominator divides a monomial times
$$\prod_{j\not\in I} {\rm Res}(f_{i_1},\dots,f_{i_d},f_j)\,.$$
Since  $\,\{ a_k \, | \, k\in I(j) \}\,$ is the unique essential subset
of $\{ a_i \, | \,i\in I\cup \{j\} \}$,
with ``essential'' as defined in \cite{rhf},  we have that 
$${\rm Res}(f_{i_1},\dots,f_{i_d},f_j) \quad = \quad
{\rm Res}(I(j);x,y).$$
We know by Proposition \ref{globalri} that $P$ is non zero.
Moreover, if any of the factors ${\rm Res}(I(j);x,y)$ were missing
from the denominator of $R_I(\one,\gamma)$, then the Laurent series
(\ref{gloresexp}) would contain only finitely many powers of $x_j$.
The formula in Proposition~\ref{globalri} implies that is impossible.  
\qed

\smallskip

For unimodular bases, 
Theorem \ref{formula2} can be refined as follows: 

\begin{proposition}\label{formula1}
Suppose that $\,\{a_i \,| \, i\in I \} \,$ is  a $\Z$-basis of $\Z^d$.    Then
\begin{equation}
R_I(\one,\gamma)\ =\ 
\frac {x^a y^b} {\prod_{j\not\in I} {\rm Res}(I(j);x,y)}
\end{equation}
where $x^a$ and $y^b$ are monomials specified in the proof.
\end{proposition}

\proof Choose $\nu , n_j\in \Z^I$, $j\in J$, so that
 $\gamma = M_I\cdot \nu$,
$a_j = M_I\cdot n_j$,  Then
$$\gamma + \sum_{j\in J} \mu_j \cdot a_j\ =\ 
M_I\cdot(\nu  + \sum_{j\in J} \mu_j \cdot m_j)
\qquad \hbox{ for all $\mu\in \N^J$},  $$
and consequently, the Laurent series 
(\ref{gloresexp}) reduces, up to sign, to
\begin{eqnarray*}
R_I(\one,\gamma) &\ =\ &
\frac {x_I^{\nu  - \one}} {y_I^{\nu }\,x_J }\ \sum_{\mu\in\N^J} 
\prod_{j\in J} \prod_{i\in I} x_i^{n_{ij}\mu_j} y_i^{-n_{ij}\mu_j}
y_j^{\mu_j} x_j^{-\mu_j}\\
&\ =\ &
\frac {x_I^{\nu  - \one}} {y_I^{\nu }}\ 
\prod_{n_{ij}>0} y_i^{n_{ij}}
\ \prod_{n_{ij}<0} x_i^{- n_{ij}}\  \prod_{j\in J} {\rm
Res}(I(j);x,y)^{-1}\ \hbox{\qed}
\end{eqnarray*}

\section{The lower bound and the linear relations}

In this section we establish the lower bound in Theorem \ref{maintheorem}
by exhibiting $| \chi(A) |$ many linearly independent binomial residues
$R_I(\beta,\gamma)$ for fixed $\beta,\gamma$ and fixed Lawrence matrix
\begin{equation*}
A\quad :=\quad
\left(
\begin{array}{cc}
I_n & I_n\\
 0 \, 0 \, \cdots \, 0
& a_1 \, a_2 \, \cdots \, a_n \\
\end{array}
\right).
\end{equation*}
We will show that all linear relations among the $R_I(\beta,\gamma)$
arise from Proposition \ref{stokes} and correspond to
 Orlik-Solomon relations  \cite[\S 3.1]{ot}.

The Gale dual to the Lawrence matrix $A$ has the form
\begin{equation}\label{lawrencegale}
B\ =\ \{b_1,\dots,b_{n},-b_1,\dots,-b_{n}\},
\end{equation}
where $B_0  = \{b_1,\dots,b_n\} \subset \Z^{n-d}$ is a
Gale dual of $\{a_1,\ldots,a_n \}$.
According to Corollary~\ref{lowb} and Lemma \ref{nbc}, the dimension of the
space of stable rational $A$-hypergeometric functions of 
degree $(-\beta,-\gamma)$ is at most the
number of  nbc-bases in $B$, which agrees with the
number of nbc-bases in $B_0$.  
The following converse will imply
Theorem \ref{maintheorem}.

\begin{theorem}
\label{dualnbcworks}
Let $\beta \in \Z_{> 0}^n$ and $\gamma \in \Z^d$.
Then the set of binomial residues $R_I(\beta,\gamma)$,
where $\{1,\ldots\! ,n\} \backslash I$ runs over all
nbc-bases of $B_0$, is linearly independent modulo
the space of unstable rational functions.
\end{theorem}

It is convenient to use the following characterization
for being an nbc-basis of the dual matroid.
The proof of Lemma \ref{replacement} is straightforward.

\begin{lemma}\label{replacement}
The set $\{1,\ldots\! ,n\} \backslash I$ is an nbc-basis of 
$B_0  = \{b_1,\dots,b_n\}$
if and only if, for each $i_0\in I$,
there exists $j_0\in \{1,\ldots\! ,n\} \backslash I \,$
such that  $j_0 > i_0$  and
$\, I\backslash \{i_0\} \cup \{j_0\}$
 is a basis of $\{a_1,\ldots,a_n\} \subset \R^d $.
\end{lemma}

\noindent{\sl Proof of Theorem~\ref{dualnbcworks}.\ }
Consider the space $\ScS(\beta,\gamma)$ of stable rational
hypergeometric functions defined in the Introduction.
The derivative $\partial_{x_i}$ induces a monomorphism
from $\ScS(\beta,\gamma)$ into $\ScS(\beta+e_i,\gamma)$, while 
$\partial_{y_i}$ induces an monomorphism into 
$\ScS(\beta+e_i,\gamma+a_i)$.
Binomial residues are mapped to binomial residues,
with the set of irreducible factors in their denominators preserved.
We may thus assume $\beta = \one$.
All linear spaces in this proof
 are understood modulo unstable rational functions.

By Theorem~\ref{formula2}, for any  basis $I$ of $\{a_1,\ldots,a_n\}$,
the denominator of $R_I(\one,\gamma)$ equals a
monomial multiplied by
\begin{equation}\label{denominator}
\prod_{j\not\in I} {\rm Res}(I(j);x,y)
\end{equation}

Let $\II_0$ denote the set of indices $I$ complementary to 
nbc-bases of $B_0$.  Let $R_{\II_0}$
denote the linear span of binomial residues
$R_I(\one,\gamma)$, $I\in \II_0$.
  Clearly, $n\not\in I$ for
any $I\in \II_0$.  Our goal is to show 
$\,\dim_\C (R_{\II_0}) \,= \,\# \,\II_0 $.

Let $K$ be a circuit of $\{a_1,\ldots,a_n\}$ which contains 
the index $n$. Define $R_{\II_0}(K)$ to be
the span of all binomial residues
$R_I(\one,\gamma)$ with $I\in \II_0$ and $I(n) = K$,
i.e., $K$ is the  unique circuit in $I \cup \{n\}$.
We may decompose
\begin{equation}\label{decomposition}
R_{\II_0} = \bigoplus_K R_{\II_0}(K)
\end{equation}
The sum in (\ref{decomposition}) is direct because
no element in $\sum_{K'\not= K} R_{\II_0}(K')$ contains
 ${\rm Res}(K;x,y)$ in its denominator, while all elements in
$R_{\II_0}(K)$ do.  

Thus, it suffices to fix $K = K_0$ and 
show that the binomial residues
$R_I(\one,\gamma)$ with $I\in \II_0$ and $I(n) = K_0$ are linearly
independent.  Let
$$ \II_1 \ =\ \{I\in \II_0 : I(n) = K_0\}\,.$$
Let $n_1$ denote the largest index which does not belong to $K_0$, 
then note that $n_1 \not\in I$ for any $I\in \II_1$.  Indeed, if
$n_1 \in I$, $I\in \II_1$, then we would not be able to replace
$a_{n_1}$ by $a_j$ with $j>n_1$ and still have a basis; this would
contradict Lemma~\ref{replacement}.  This means that we can repeat
the previous argument with $\II_1$ in place of $\II_0$ and
$n_1$ in place of $n$ and obtain a decomposition
of $R_{\II_1}$ as a direct sum of subspaces $R_{\II_1}(K)$ spanned
by binomial residues $R_I(\one,\gamma)$ with
$I\in \II_1$ and $I(n_1) = K$.  Continuing in this manner, all
 subspaces $R_{\II_p}(K)$ will  eventually be
one-dimensional. Then, the desired result
 follows from Proposition
\ref{globalri}.
\qed

\smallskip

We next describe all linear relations among the binomial
residues $\,R_I(\beta,\gamma) \,$ as $I$ varies.
In the identity below, it is  essential  to keep track of signs. 
Namely, if
 $\,I' \,$ is taken to be ordered then we must
multiply $\, R_{I' \cup \ell}(\beta,\gamma) \,$ by the
sign of the  permutation which orders $\,I' \cup \{\ell \}$.

\begin{theorem}
\label{orlikk}
Let  $I'$ be a $(d-1)$-subset of $ \{1,\ldots,n\}\,$
and ${\rm ind}\,I'$ the set of indices $\ell$ such that
$\{ a_\ell\} \cup \{ a_i : i \in I'\}$ is a basis of $\R^d$.
Then
$$ \sum_{\ell \in {\rm ind} {I'}} R_{I' \cup \ell}(\beta,\gamma) 
\,\,\equiv \,\, 0 \qquad \hbox{modulo unstable rational functions} ,$$
and these span all the $\C$-linear relations
 relations among the $ R_{I}(\beta,\gamma) $.
\end{theorem}

\proof
By Proposition~\ref{globalri}, 
all $\,R_I(\beta,\gamma)\,$ residues are stable.
We have  established that  the spaces 
$\ScS(\beta,\gamma)$ have the same dimension $|\chi(A)|$ for all
$\beta$ and $\gamma$. It follows that the maps
$\partial_{x_i} : \ScS(\beta,\gamma) \to 
\ScS(\beta+ e_i,\gamma) $ and
$\partial_{y_i} : \ScS(\beta,\gamma) \to \ScS(\beta+ e_i,\gamma+a_i) $
are isomorphisms.  Iterating, we can   assume that $(-\beta,
-\gamma)
$ lies in the Euler-Jacobi cone $-{\rm Int}({\rm pos}(A))$.  
By Proposition \ref{ejprop}, 
there are no unstable rational $A$-hypergeometric functions,
so we are claiming that
$\, \sum_{\ell \in {\rm ind} {I'}} R_{I' \cup \ell}(\beta,\gamma) \,$
is zero.

We may assume that $\, \{ a_i : i \in I'\}\,$ is linearly independent.
On the $B$-side, the complement of ${I'}$ has $n-d+1$ elements and 
therefore defines a dependent set  $\{b_i , i \not\in I'\}$. 
We can consider as in \S 2, the central hyperplane arrangement $\scA$ 
defined by $\HH$. Consider the socle of the {\it Orlik-Solomon algebra}
of that hyperplane arrangement \cite[\S 3.1]{ot}. The 
linear relation in Theorem \ref{orlikk} is the translation
to the $A$-side of the relation in the socle degree of the
Orlik-Solomon algebra defined by   $\{b_i , i \not\in I'\}$. 
In view of  \cite[Theorem 3.4]{ot}  and Theorem \ref{maintheorem}, 
it suffices to show that the asserted relations are valid. It will then follow
by dimension reasons that they span all $\C$-linear relations.

We now prove the identity 
$\, \sum_{\ell \in {\rm ind} {I'}} R_{I' \cup \ell}(\beta,\gamma) \, = \, 0\,$
using the formulation in terms of  toric residues 
given in \S 2. By (\ref{eje}),
all $h_j (\beta,\gamma)$ are non negative, and so the polar divisor of the
form $\Phi(\beta,\gamma)$ in (\ref{meromtoric}) is contained in the union of 
the divisors $Y_i = \{F_i =0 \}, i=1,\dots,n.$

For $k= 1,\dots,d-1$, set $G^{I'}_k = F_{i_k}$. Set also $G^{I'}_d =
\prod_{j \notin I'} F_j$ and let $G^{I'}_0 = z_1 \dots z_{2p}.$  Then,
$G^{I'}_0,\dots, G^{I'}_d$ define divisors with empty intersection in
$X =X_\Delta$ for generic values of the coefficients and moreover 
$$\Phi(\beta,\gamma)  \quad = \quad  
\frac {z^{h(\beta,\gamma)}\ \ \Omega_\Delta} {G_1 \dots
G_d}.$$
Proposition \ref{stokes} implies that the  corresponding toric residue
vanishes:
$${\rm Res}^X_{G^{I'}}(\Phi(\beta,\gamma)) \quad = \quad 0.$$
On the other hand,  consider also the following $n-d +1$ families
of divisors: for any $\ell \notin {I'}$, set $G^{{I'},\ell}_k = G^{I'}_k$ for any $k = 1,\dots, d-1,$
$G^{{I'},\ell}_d= F_\ell$ and $G^{{I'},\ell}_0 =
\prod_{j \notin {I'} \cup \{\ell\}} F_j.$  Again, these divisors have empty intersection
on $X$ for generic values of the coefficients
and the poles of $\Phi(\beta,\gamma)$ are contained in their union, and so
we can consider the toric residues ${\rm Res}^X_{G^{{I'},\ell}}(\Phi(\beta,\gamma)).$
These toric residues are non-zero precisely when 
$\ell \in {\rm ind} \,I'$. 
We conclude that the  following relations hold:
\begin{eqnarray*}
\sum_{\ell \in {\rm ind} {I'}}
\! {\rm Res}^X_{G^{{I'},\ell}}(\Phi(\beta,\gamma))
\, = \,
\sum_{\ell \notin j} {\rm Res}^X_{G^{{I'},\ell}}(\Phi(\beta,\gamma))
\, = \, {\rm Res}^X_{G^{I'}}(\Phi(\beta,\gamma))  \, = \,0.
 \end{eqnarray*}
The second equality follows from a variation on \cite[\S II.7]{tsikh}.
Translating back to binomial residues
 completes the proof of Theorem \ref{orlikk}.
\ \qed

\medskip
In \cite{rhf}, we studied the
problem of classifying vector configurations $A$ 
for which there exist rational a $A$-hypergeometric function
which is not a Laurent polynomial.
We conjectured \cite[Conjecture
1.3]{rhf} that such a configuration has to have a facial subset which
is an essential Cayley configuration.  It is easy to see that
Lawrence liftings  are  Cayley configurations of segments; they are
essential if and only if $n = d+1$.
We also conjectured  \cite[Conjecture
5.7]{rhf} that a rational $A$-hypergeometric function has
an iterated derivative which is a linear combination of toric residues
associated with facial subsets of $A$.  

\begin{theorem}
Conjecture~5.7 in \cite{rhf} holds for Lawrence configurations.
\end{theorem}

\proof
Let $A$ be a Lawrence configuration.  The assertion of
\cite[Conjecture~5.7]{rhf} is obvious for unstable rational
hypergeometric functions.  On the other hand, given a stable
rational hypergeometric function, a suitable derivative will
have degree in the Euler-Jacobi cone and hence, by
Theorem~\ref{maintheorem}, will be a linear combination of toric
residues. \qed

\medskip
\medskip

\noindent{\bf Acknowledgements:}  
 Alicia Dickenstein was partially
supported by UBACYT TX94 and CONICET, Argentina,
and the Wenner-Gren Foundation, Sweden.
Bernd Sturmfels was partially supported
by NSF Grant DMS-9970254.

\medskip 
\medskip 
\medskip 
\medskip 

\end{document}